\documentclass[12pt]{amsart}
\usepackage{amssymb}

\newcommand{\Bgp}{{\Z^\N}}
\newcommand{\Arh}{Arhangel'ski\u{\i}}

\long\def\forget#1\forgotten{}
\newcommand{\issuenumber}{22}
\newcommand{\issuemonth}{September}
\newcommand{\issueyear}{2007}

\newtheorem{thm}{Theorem}[section]
\newtheorem{prob}[thm]{Problem}

\newtheorem{issue}{Issue}

\theoremstyle{definition}
\newtheorem{defn}[thm]{Definition}

\theoremstyle{remark}

\newcommand{\ed}{
\newpage

\section{Unsolved problems from earlier issues}

\stepcounter{issue}

\stepcounter{issue}

\stepcounter{issue}

\begin{issue}
Does $\sone(\Omega,\Tau)$ imply $\ufin(\Gamma,\Gamma)$?
\end{issue}

\begin{issue}
Is $\fp=\fp^*$? (See the definition of $\fp^*$ in that issue.)
\end{issue}

\begin{issue}
Does there exist (in ZFC) an uncountable set satisfying $\sone(\BG,\B)$?
\end{issue}

\stepcounter{issue}

\begin{issue}
Does $X \nin \NON(\M)$ and $Y\nin\mathsf{D}$ imply that
$X\cup Y\nin \COF(\M)$?
\end{issue}

\begin{issue}[CH]
Is $\split(\Lambda,\Lambda)$ preserved under finite unions?
\end{issue}

\begin{issue}
Is $\cov(\M)=\fo$? (See the definition of $\fo$ in that issue.)
\end{issue}

\begin{issue}
Does $\sone(\Gamma,\Gamma)$ always contain an element of cardinality $\fb$?
\end{issue}

\begin{issue}
Could there be a Baire metric space $M$ of weight $\aleph_1$ and a partition
$\mathcal{U}$ of $M$ into $\aleph_1$ meager sets where for each ${\mathcal U}'\subset\mathcal U$,
$\bigcup {\mathcal U}'$ has the Baire property in $M$?
\end{issue}

\stepcounter{issue} 

\begin{issue}
Does there exist (in ZFC) a set of reals $X$ of cardinality $\fd$ such that all
finite powers of $X$ have Menger's property $\sfin(\cO,\cO)$?
\end{issue}

\begin{issue}
Can a Borel non-$\sigma$-compact group be generated by a Hurewicz subspace?
\end{issue}

\begin{issue}[MA]
Is there an uncountable $X\sbst\R$ satisfying $\sone(\BO,\BG)$?
\end{issue}

\begin{issue}[CH]
Is there a totally imperfect $X$ satisfying $\ufin(\cO,\Gamma)$
that can be mapped continuously onto $\Cantor$?
\end{issue}

\begin{issue}[CH]
Is there a Hurewicz $X$ such that $X^2$ is Menger but not Hurewicz?
\end{issue}

\begin{issue}
Does the Pytkeev property of $C_p(X)$ imply the Menger property of $X$?
\end{issue}

\begin{issue}
Does every hereditarily Hurewicz space satisfy $\sone(\BG,\BG)$?
\end{issue}

\begin{issue}[CH]
Is there a Rothberger-bounded $G\le\Bgp$ such that $G^2$ is not Menger-bounded?
\end{issue}

\begin{issue}
Let $\cW$ be the van der Waerden ideal.
Are $\cW$-ultrafilters closed under products?
\end{issue}

\general\end{document}}

\newcommand{\Cantor}{{\{0,1\}^\N}}

\newcommand{\fb}{\mathfrak{b}}
\newcommand{\fg}{\mathfrak{g}}

\newcommand{\fd}{\mathfrak{d}}
\newcommand{\fp}{\mathfrak{p}}

\newcommand{\NON}{{\mathsf   {NON}}}
\newcommand{\COF}{{\mathsf   {COF}}}

\newcommand{\cI}{\mathcal{I}}

\newcommand{\M}{\mathcal{M}}

\newcommand{\cov}{\mathsf{cov}}

\newcommand{\R}{\mathbb{R}}

\newcommand{\fo}{\mathfrak{od}}

\renewcommand{\split}{\mathsf{Split}}
\newcommand{\bq}{\begin{quote}}
\newcommand{\eq}{\end{quote}}
\newcommand{\cO}{\mathcal{O}}
\newcommand{\B}{\mathcal{B}}
\newcommand{\BG}{\B_\Gamma}

\newcommand{\BO}{\B_\Omega}

\newcommand{\sone}{\mathsf{S}_1}    \newcommand{\sfin}{\mathsf{S}_{fin}}

\newcommand{\ufin}{\mathsf{U}_{fin}}

\newcommand{\nin}{\not\in}

\newcommand{\cU}{\mathcal{U}}
\newcommand{\cV}{\mathcal{V}}
\newcommand{\cW}{\mathcal{W}}

\newcommand{\N}{\mathbb{N}}
\newcommand{\Z}{\mathbb{Z}}

\newcommand{\sbst}{\subseteq}
\newcommand{\by}[2]{\par\hfill\emph{#1}, #2}
\newcommand{\nby}[1]{\par\hfill\emph{#1}}
\newcommand{\Tau}{\mathrm{T}}
\newcommand{\CE}{\textsc{CE}}

\newcommand{\be}{\begin{enumerate}}
\newcommand{\ee}{\end{enumerate}}
\newcommand{\bi}{\begin{itemize}}
\newcommand{\ei}{\end{itemize}}

\newcommand{\general}{\small\vfill\par\noindent\hrulefill\par
\noindent\textbf{Previous issues.} The previous issues of this
bulletin, which contain general information (first issue), basic
definitions, research announcements, and open problems (all
issues) are available online,
at \texttt{http://front.math.ucdavis.edu/search?\&t=\%22SPM+Bulletin\%22}
\\[0.1cm]
\textbf{Contributions.}
Please submit your contributions (announcements, discussions, and open problems)
by e-mailing us. It is preferred to write them
in \LaTeX{}.
The authors are urged to use as standard notation as possible, or otherwise give
the definitions or a reference to where the notation is explained.
Contributions to this bulletin would not require any transfer of copyright,
and material presented here can be published elsewhere.\\[0.1cm]
\textbf{Subscription.}
To receive this bulletin (free) to your
e-mailbox, e-mail us:\\
{boaz.tsaban@weizmann.ac.il}
}

\newcommand{\link}[1]{\par\hfill{\texttt{#1}}}

\newcommand{\nArxPaper}[5]{\subsection{#2}{#4}\par\hfill{\arx{#1}}\par\hfill\emph{#3}}

\newcommand{\arx}[1]{\texttt{http://arxiv.org/abs/#1}}
\newcommand{\url}[1]{\bq\texttt{#1}\eq}
\newcommand{\online}[1]{The paper is available online at \url{#1}}

\title[$\mathcal{SPM}$ Bulletin \textbf{\issuenumber} (\issuemonth{} \issueyear)]{%
$\mathcal{SPM}$ Bulletin\\[0.5cm]
Issue number \issuenumber: \issuemonth{} \issueyear{} \CE{}}

\begin{document}
\maketitle

\tableofcontents

\section{Editor's note}

We are now after the \emph{First European Set Theory Meeting},
a historically important and well organized event. The talks gave
the right blend of theory and applications of set theory.
Thanks to Benedikt Loewe, Grzegorz Plebanek, Jouko V\"a\"an\"anen, and Boban Velickovic
for the organization.

Following Jana Fla\v{s}kova's talk at this meeting, we have invited her to
contribute a section to this issue.
We thank her for her interesting contribution in Second \ref{flas} and
in the \emph{Problem of the Issue} section.

The list of problems at the end of the bulletin became longer than one page.
We therefore removed the first few, and will continue this way unless some more
problems are solved and their space becomes available\dots

A much better version of Shelah's paper showing that $\fg\le\fb^+$ is now available
at \texttt{arxiv.org/abs/math/0612353}

\medskip

Enjoy,

\by{Boaz Tsaban}{boaz.tsaban@weizmann.ac.il}

\hfill \texttt{http://www.cs.biu.ac.il/\~{}tsaban}

\section{Invited contribution: Ultrafilters and small sets}\label{flas}

There have been several attempts to connect ultrafilters with
families of ``small" sets. Two of them --- $0$-points and
$\mathcal{I}$-ultrafilters --- were important for my Ph.D. thesis.
The first one is due to Gryzlov \cite{Gr}: an ultrafilter
$\mathcal{U} \in \N^{\ast}$ is called a \emph{$0$-point} if for
every one-to-one function $f: \N \rightarrow \N$ there exists a
set $U \in \mathcal{U}$ such that $f[U]$ has asymptotic density
zero. The second term was introduced by Baumgartner \cite{B}: Let
$\mathcal{I}$ be a~family of subsets of a~set $X$ such that
$\mathcal{I}$ contains all singletons and is closed under subsets.
Given a free ultrafilter $\mathcal{U}$ on $\N$, we say that
${\mathcal{U}}$ is an $\mathcal{I}$-ultrafilter if for any
$F:\N\rightarrow X$ there is $A \in \mathcal{U}$ such that
$F[A] \in \mathcal{I}$.

In my Ph.D. thesis \cite{FTh} (on which all my papers are more
less based) I studied $\mathcal{I}$-ultrafilters in the setting $X
= \N$ and $\mathcal{I}$ is an ideal on $\N$ or another
family of ``small" subsets of $\N$ that contains finite sets and
is closed under subsets. As $\mathcal{I}$ were considered the
ideal of sets with asymptotic density zero
{${\mathcal{Z}_0}$}$=\{A \subseteq \N: \limsup_{n \rightarrow
\infty} \frac{|A \cap n|}{n} = 0\}$, the summable ideal
{${\mathcal{I}_{1/n}}$}$= \{A \subseteq \N: \sum_{a \in A}
\frac{1}{a} < \infty\}$ or the family of (almost) thin sets and
$(SC)$-sets.

Here are the (probably not common) definitions: We say that $A
\subseteq \N$ with an increasing enumeration $A = \{a_n:n \in
\N\}$ is
\begin{description}
\item[\emph{thin}] if $\lim_{n \rightarrow \infty} {\frac{a_n}{a_{n+1}}} = 0$;
\item[\emph{almost thin}] if $\lim_{n \rightarrow \infty} {\frac{a_n}{a_{n+1}}} < 1$;
\item [\emph{$(SC)$-set}] if $\lim_{n \rightarrow \infty} a_{n+1}-a_n = \infty$.
\end{description}

In the thesis various examples of $\mathcal{I}$-ultrafilters for
all these (and also some other) families $\mathcal{I}$ are
constructed under additional set theoretic assumptions.

In my first paper \cite{F1} it is shown that thin sets and almost
thin sets actually determine the same class of
$\mathcal{I}$-ultrafilters and there is a proof that the existence
of these ultrafilters is independent of ZFC. The relation between
this class of ultrafilters and selective ultrafilters or
$Q$-points is studied. Some construction made in the paper under
CH were proved in the thesis assuming MA${}_{\rm{ctble}}$.

The next paper \cite{F2} focuses on $\mathcal{I}$-ultrafilters
where $\mathcal{I}$ is the summable ideal ${\mathcal{I}_{1/n}}$ or
the density ideal ${\mathcal{Z}_0}$. The relation between these
two classes of ultrafilters is shown and also the relation to the
class of $P$-points. Assuming CH or MA${}_{\rm{ctble}}$ several
examples of these ultrafilters are constructed. Again, stronger
versions of some of the results can be found in the thesis.

One of the few ZFC results in my thesis is the following: There
exists an ultrafilter $\mathcal{U} \in \N^{\ast}$ such that for
every one-to-one function $f: \N \rightarrow \N$ there exists a
set $U \in \mathcal{U}$ with $f[U]$ in the summable ideal. This
theorem strengthens Gryzlov's result concerning the existence of
$0$-points and it was published also separately in \cite{F3}.

Connections between various $\mathcal{I}$-ultrafilters and some
well-known ultrafilters such as $P$-points were studied in two
sections of my thesis. It is known that $P$-points can be
described as $\mathcal{I}$-ultrafilters in two different ways: If
$X = 2^{\N}$ then $P$-points are precisely the
$\mathcal{I}$-ultrafilters for $\mathcal{I}$ consisting of all
finite and converging sequences, if $X = \omega_1$ then $P$-points
are precisely the $\mathcal{I}$-ultrafilters for $\mathcal{I} =
\{A \subseteq \omega_1: A \hbox{ has order type } \leq \omega\}$.
My latest paper \cite{F4} deals with the question whether there is
a family $\mathcal{I}$ of subsets of natural numbers such that
$P$-points are precisely the $\mathcal{I}$-ultrafilters. However,
only some partial answers are given.

During the 1st European Set Theory Meeting in B\c{e}dlewo I gave a
talk ``On sums and products of certain
$\mathcal{I}$-ultrafilters". As the title suggests it was a
summary of my knowledge about sums and products of some
$\mathcal{I}$-ultrafilters. The slides and notes with proofs on
which the talk was based (as well as my Ph.D. thesis) are
available online on my webpage:

{\tt http://home.zcu.cz/\~{}flaskova}

\nby{Jana Fla\v{s}kov\'a}

\hfill Department of Mathematics, University of West Bohemia

\section{Research announcements}

\nArxPaper{0706.3815}
{Inverse Systems and I-Favorable Spaces}
{Andrzej Kucharski and Szymon Plewik}
{Let $X$ be a compact space.
Player I has a winning strategy in the open-open game played on $X$
if, and only if $X$ can be represented as a
limit of $\sigma$-complete inverse system of compact metrizable spaces with
skeletal bonding maps.}

\nArxPaper{0706.3729}
{Combinatorial and hybrid principles for $\sigma$-directed families of countable sets modulo finite}
{James Hirschorn}
{We consider strong combinatorial principles for $\sigma$-directed families of
countable sets in the ordering by inclusion modulo finite, e.g. $P$-ideals of
countable sets. We try for principles as strong as possible while remaining
compatible with CH, and we also consider principles compatible with the
existence of nonspecial Aronszajn trees. The main thrust is towards abstract
principles with game theoretic formulations. Some of these principles are
purely combinatorial, while the ultimate principles are primarily combinatorial
but also have aspects of forcing axioms.
}

\nArxPaper{0707.1313}
{A dichotomy characterizing analytic digraphs of uncountable Borel chromatic number in any dimension}
{Dominique Lecomte}
{We study the extension of the Kechris-Solecki-Todorcevic dichotomy on
analytic graphs to dimensions higher than $2$. We prove that the
extension is possible in any dimension, finite or infinite. The
original proof works in the case of the finite dimension. We first
prove that the natural extension does not work in the case of the
infinite dimension, for the notion of continuous homomorphism used
in the original theorem. Then we solve the problem in the case of
the infinite dimension. Finally, we prove that the natural
extension works in the case of the infinite dimension, but for the
notion of Baire-measurable homomorphism.}

\nArxPaper{0707.1313}
{A dichotomy characterizing analytic digraphs of uncountable Borel chromatic number in any dimension}
{Dominique Lecomte}
{We study the extension of the Kechris-Solecki-Todorcevic dichotomy on
analytic graphs to dimensions higher than $2$. We prove that the extension is
possible in any dimension, finite or infinite. The original proof works in the
case of the finite dimension. We first prove that the natural extension does
not work in the case of the infinite dimension, for the notion of continuous
homomorphism used in the original theorem. Then we solve the problem in the
case of the infinite dimension. Finally, we prove that the natural extension
works in the case of the infinite dimension, but for the notion of
Baire-measurable homomorphism.}

\nArxPaper{0707.1818}
{Large continuum, oracles}
{Saharon Shelah}
{Our main theorem is about iterated forcing for making the continuum larger
than $\aleph_2$. We present a generalization of \texttt{math.LO/03032\allowbreak 94} which is dealing
with oracles for random, etc., replacing $\aleph_1,\aleph_2$ by $\lambda,\lambda^+$
(starting with $\lambda=\lambda^{<\lambda}>\aleph_1$). Well, instead of properness we
demand absolute c.c.c. So we get, e.g. the continuum is $\lambda^+$ but we can get
$\cov(\M)=\lambda$. We give some applications. As in \texttt{math.LO/0303294}, it is a
``partial'' countable support iteration but it is c.c.c.}

\nArxPaper{0707.1967}
{Borel hierarchies in infinite products of Polish spaces}
{Rana Barua and Ashok Maitra}
{Let $H$ be a product of countably infinite number of copies of an uncountable
Polish space $X$. Let $\Sigma_\xi$ $(\bar {\Sigma}_\xi)$ be the class of Borel
sets of additive class $\xi$ for the product of copies of the discrete topology
on $X$ (the Polish topology on $X$), and let ${\mathcal B} = \cup_{\xi < \omega_1}
\bar{\Sigma}_\xi$. We prove in the L\'{e}vy--Solovay model that
$\bar{\Sigma}_\xi =\Sigma_{\xi}\cap {\mathcal B}$ for $1 \leq \xi < \omega_1$.}

\subsection{A game for the Borel functions}
Abstract. We present an infinite game that characterizes the
Borel functions on Baire Space.
\link{www.illc.uva.nl/Publications/ResearchReports/PP-2006-24.text.pdf}
\nby{Brian Semmes}

\nArxPaper{0708.1981}
{On some problems in general topology}
{Saharon Shelah}
{We prove that \Arh's problem has a consistent positive answer: If
$V$ is a model of CH, then for some $\aleph_1$-complete $\aleph_2$-c.c.\ forcing notion $P$
of cardinality $\aleph_2$, we have that $P$ forces ``CH and there is a Lindel\"of regular
topological space of size $\aleph_2$ with clopen basis with every point of
pseudo-character $\aleph_0$ (i.e., each singleton is the intersection of countably
many open sets)''.

Also, we prove the consistency of: CH + $2^{\aleph_1} > \aleph_2$ + ``there is no
space as above with $\aleph_2$ points'' (starting with a weakly compact cardinal).

Appeared in: \emph{Set Theory}, Boise ID, 1992--1994, Contemporary Mathematics, vol.\ 192, 91--101.}

\section{Problem of the Issue}

Definitions not stated below can be found in Section \ref{flas} above.
Some more definitions concerning ultrafilters can be found in \cite{CN}.

The problem of this issue concerns products of $\cI$-ultrafilters
for the case $X = \N$ and $\cI$ is an ideal on $\N$.

\begin{defn}
If $\cU$ and $\cV$ are ultrafilters on $\N$ then $\cU \cdot \cV$ is
the ultrafilter on $\N \times \N$ defined by $M \in \cU \cdot \cV$
if and only if $\{n : \{m: \langle n, m \rangle \in M\} \in \cV\}
\in \cU$. Since isomorphic ultrafilters can be identified we may
regard $\cU \cdot \cV$ as an ultrafilter on $\N$. The ultrafilter
$\cU \cdot \cV$ is called the \emph{product of ultrafilters $\cU$
and} $\cV$.
\end{defn}

\begin{defn} Let $\cI$ be an ideal on $\N$. We say that $\cI$-ultrafilters are \emph{closed
under products} if the product of two arbitrary $\cI$-ultrafilters
is again an $\cI$-ultrafilter.
\end{defn}

For a $P$-ideal $\cI$ the class of $\cI$-ultrafilters is closed
under products \cite{FTh}. However, not much is known for other
ideals. For example, if $\cI$ is the ideal generated by thin sets or
the ideal generated by $(SC)$-sets then $\cI$-ultrafilters are not
closed under products \cite{FTh}. In fact, even more is true.

\begin{thm}[\cite{FTh}]
For every $\cU \in \N^{\ast}$ the ultrafilter $\cU \cdot \cU$ is
not an $(SC)$-ultrafilter (thin ultrafilter).
\end{thm}

This property shares the class of all $P$-points (the partition
$\{\{n\} \times \N: n \in \N\}$ witnesses the fact that no product
of two free ultrafilters is a $P$-point), but it is consistent with
ZFC that there exist thin ultrafilters (and hence
$(SC)$-ultrafilters) which are not $P$-points \cite{FTh}.

Another example of an ideal which is not a $P$-ideal is the following.

\begin{defn}
The \emph{van der Waerden ideal} $\cW$ is the family of all $A
\subseteq \N$ such that $A$ does \emph{not} contain arithmetic
progressions of arbitrary length.
\end{defn}

\begin{prob}
Are $\cW$-ultrafilters closed under products?
\end{prob}

A positive answer would provide (consistent) examples of
$\cW$-ultrafilters that are neither $P$-points nor (SC)-ultrafilters
(and thin ultrafilters). But I expect rather a negative answer.

\nby{Jana Fla\v{s}kov\'a}

\ed